\documentclass[11pt,a4paper]{amsart}

\usepackage[T1]{fontenc}
\usepackage{lmodern}
\usepackage{amsmath,amssymb,amsthm,mathrsfs}
\usepackage[margin=1in]{geometry}
\usepackage[colorlinks=true,linkcolor=blue,citecolor=blue,urlcolor=blue]{hyperref}

\newtheorem{theorem}{Theorem}[section]
\newtheorem{proposition}[theorem]{Proposition}
\newtheorem{lemma}[theorem]{Lemma}

\theoremstyle{definition}
\newtheorem{definition}[theorem]{Definition}
\theoremstyle{remark}
\newtheorem{remark}[theorem]{Remark}

\numberwithin{equation}{section}
\newcommand{\cM}{\mathcal M}
\newcommand{\rk}{\operatorname{rk}}
\newcommand{\Alb}{\operatorname{Alb}}
\newcommand{\Hom}{\operatorname{Hom}}
\newcommand{\AJ}{\operatorname{AJ}}

\title[Rank of the Fano normal function]{The rank of the Fano normal function via the Fermat cubic}
\author[L. Fassina]{Lorenzo Fassina}

\address{Dipartimento di Matematica, Universit\`a degli Studi di Pavia, Via Ferrata 5, 27100 Pavia, Italy}

\email{lorenzo.fassina02@universitadipavia.it}
\date{}



\begin{document}

\begin{abstract}
We prove that the canonical derivative of the Fano normal
function has rank ten at the Fermat cubic threefold.
This yields a direct proof of its generic maximal rank, a result
also implied by a general theorem of Gao--Zhang.
Our approach combines Hain's equivariance principle with
an explicit representation-theoretic computation and the
nonvanishing of the Griffiths infinitesimal invariant established
by Collino--Naranjo--Pirola.
\end{abstract}

\maketitle

% -------------------------------------------------
\section{Introduction}
\label{sec:introduction}
% -------------------------------------------------

Let $V\subset\mathbb P^4$ be a smooth cubic threefold and let
$F(V)$ denote its Fano surface of lines. The 
isomorphism $\Alb(F(V))\simeq J(V)$ allows us, after
choosing a base point, to regard $F(V)$ as a surface
in the intermediate Jacobian $J(V)$. Its image under
the involution $x\mapsto-x$ will be denoted by $F(V)^-$.
The difference
\[
    Z_V:=F(V)-F(V)^-
\]
is a homologically trivial algebraic cycle on $J(V)$.
As shown by Collino--Naranjo--Pirola \cite{CNP}, the
primitive component of its Abel--Jacobi class is
independent of the chosen base point. Varying $V$ in
the moduli space $\cM$ of smooth cubic threefolds,
these classes define the \emph{Fano normal function}
$\nu_F$.

Normal functions arise naturally from families of
homologically trivial algebraic cycles through the
Abel--Jacobi map. Their infinitesimal behaviour
was studied by Griffiths \cite{Griffiths1983},
who introduced the infinitesimal invariant,
and subsequently by Green \cite{Green1989}
and Voisin \cite{Voisin1988}.
For cubic threefolds, Collino--Naranjo--Pirola
\cite{CNP} used this invariant to study the
Fano normal function, proving the non-torsion
of its primitive component and showing that
its infinitesimal invariant determines the cubic.

Another important aspect of normal functions is
their rank, which measures their variation with
respect to the locally constant real structure
of the family of intermediate Jacobians.
This is closely related to the study of Betti maps
for sections of abelian schemes.
Andr\'e--Corvaja--Zannier \cite{ACZ}
related the rank of the Betti map to the
Kodaira--Spencer map, obtaining results on
its maximality and applications to torsion points.
Further results on the generic rank of Betti maps
were obtained by Gao and Zhang \cite{GaoBetti}, \cite{GaoZhang2024}.

In addition, Hain \cite{Hain} independently established the maximal rank of the Ceresa normal function over the moduli space of curves of genus
$g\geq 3$, using different methods.
These developments motivate the study of
the rank of the Fano normal function
associated with cubic threefolds.

In this paper we give a direct computation of
the rank of the Fano normal function.
Our main result is the following.

\begin{theorem}\label{thm:main}
The Fano normal function has maximal rank:
\[
    \operatorname{rk}(\nu_F)=10=\dim\cM.
\]
More precisely, its canonical derivative has rank $10$
at the Fermat cubic threefold.
\end{theorem}

The generic maximal rank of the Fano normal function also follows
from a general theorem of Gao--Zhang on the Betti rank of normal
functions. We recall their theorem and explain why it applies
in our setting.

\begin{theorem}[Gao--Zhang {\cite[Theorem~5.1]{GaoZhang}}]
\label{thm:Gao-Zhang}

Let $S$ be a smooth quasi-projective variety and
let $\nu$ be an admissible normal function associated with a polarized variation of Hodge structure $\mathbb V_{\mathbb Q}$ of weight $-1$. Assume that
$\mathbb V_{\mathbb Q}$ is irreducible and that $\nu$ is not a torsion
section. Then
\[
    \rk(\nu)
    =
    \min\left\{
        \dim\Phi_\nu(S),\,
        \frac12\dim_{\mathbb Q}\mathbb V_{\mathbb Q,s}
    \right\},
\]
where $\Phi_\nu$ denotes the period map of the mixed variation of
Hodge structure associated with $\nu$.
\end{theorem}

Let us explain why Gao--Zhang's theorem applies to the Fano
normal function. Its underlying rational variation is $\mathbb V_V=P^5(J(V),\mathbb Q)(3),$
which has weight $-1$. The variation $\mathbb V$ is irreducible. Indeed, the monodromy of the family of cubic threefolds is Zariski-dense in
$\operatorname{Sp}_{10}$ by \cite{Beauville},
and the primitive fifth exterior power of the
standard symplectic representation is irreducible.
Moreover, $\operatorname{rank}_{\mathbb Q}\mathbb V
    =\binom{10}{5}-\binom{10}{3}=132.$
The Fano normal function is admissible, being induced by
algebraic cycles, and its primitive component is non-torsion
by \cite[Theorem~5.11]{CNP}.
Finally, the Torelli theorem for cubic threefolds
\cite{ClemensGriffiths} implies that the period map
$[V]\mapsto[J(V)]$ is an injective immersion and the same holds infinitesimally for the primitive fifth
cohomology.
% The same holds infinitesimally for the primitive fifth
% cohomology: indeed, the variation of the line
% \[
%     P^{5,0}(J(V))
%     =\bigwedge\nolimits^5 H^{1,0}(J(V))
% \]
% detects the infinitesimal variation of $H^{1,0}(J(V))$.
Thus the pure period map of $\mathbb V$ has image of
dimension $10$. Since the mixed period map
$\Phi_{\nu_F}$ projects to the pure period map, we obtain $\dim\Phi_{\nu_F}(\mathcal M)=10$.
Gao--Zhang's theorem therefore gives
\[
    \operatorname{rk}(\nu_F)=\min\{10,66\}=10.
\]

Our purpose is to obtain this conclusion by a different, explicit
argument, which also detects maximal rank at the special Fermat
fiber. We follow Hain's method for the Ceresa normal function
\cite{Hain}: the rank is computed by the canonical derivative, whose
construction from the real flat trivialization makes it compatible
with automorphisms. At a fixed point of a group action, the canonical
derivative is an equivariant linear map.

We apply this observation to the Fermat cubic
\[
    V_F=\{x_0^3+x_1^3+x_2^3+x_3^3+x_4^3=0\}
    \subset\mathbb P^4
\]
and its symmetry group $G\simeq(\mu_3)^4\rtimes S_5$. Let $0\in Y$
denote the Fermat point in a $G$-equivariant local deformation space.
Writing $P^{2,3}(J(V_F))$ for the primitive Hodge component of type
$(2,3)$, we use the Jacobian ring to construct explicitly a generator
\[
    \varphi\in\Hom_G\bigl(T_0Y,P^{2,3}(J(V_F))\bigr)
\]
and prove that
\[
    \Hom_G\bigl(T_0Y,P^{2,3}(J(V_F))\bigr)
    =\mathbb C\varphi,
    \qquad \rk(\varphi)=10.
\]
Using the character decomposition and the action
of the stabilizer, we construct an explicit
$G$-equivariant homomorphism and show that it
generates the entire space and has rank ten.
Since the canonical derivative
of $\nu_F$ at $V_F$ belongs to this space, maximal rank follows once
its non-vanishing is established. 
For this final step, we apply the nonvanishing
criterion of Collino--Naranjo--Pirola
\cite[Proposition~5.9(c)]{CNP} to a special
line in $\mathbb P^4$ meeting the Fermat cubic
in three distinct points.

The paper is organized as follows. In Section~2, we recall
the canonical derivative of a normal function and its basic
properties. Section~3 is devoted to the construction of the
primitive Fano normal function. In Section~4, we study the
symmetries of the Fermat cubic and determine the space of
equivariant homomorphisms. Finally, in Section~5, we prove
the non-vanishing of the canonical derivative and establish
our main result.

\section{Normal functions, rank and symmetries}
\label{sec:normal-functions}

We recall the basic constructions concerning normal functions
and their infinitesimal invariants, following
\cite[Chapter~7]{VoisinII}. For the notions of rank and canonical
derivative, we refer to \cite[Sections~2--3]{Hain}.
Let $\pi\colon\mathcal X\to Y$ be a smooth projective morphism over
a complex manifold $Y$, and write $X_y:=\pi^{-1}(y)$. Fix $p\geq 1$
and set
\[
    \mathcal H
    :=R^{2p-1}\pi_*\mathbb C\otimes_{\mathbb C}\mathcal O_Y,
    \qquad
    \mathbb H_{\mathbb Z}
    :=R^{2p-1}\pi_*\mathbb Z/\mathrm{torsion}.
\]
Thus $\mathcal H$ is a holomorphic vector bundle with fiber
$\mathcal H_y=H^{2p-1}(X_y,\mathbb C)$, while $\mathbb H_{\mathbb Z}$
is its integral local system modulo torsion. Let
\[
    \nabla\colon\mathcal H\longrightarrow
    \mathcal H\otimes\Omega_Y^1
\]
be the flat Gauss--Manin connection on the associated holomorphic bundle $\mathcal H$. The Hodge filtrations of the fibers define holomorphic
subbundles $F^i\mathcal H\subset\mathcal H$, with
\[
    F^i\mathcal H_y
    =\bigoplus_{a\geq i}H^{a,2p-1-a}(X_y).
\]
These satisfy Griffiths transversality
\cite[\S 5.1.2]{VoisinII}:
\begin{equation}\label{eq:Griffiths-transversality}
    \nabla(F^i\mathcal H)
    \subset
    F^{i-1}\mathcal H\otimes\Omega_Y^1.
\end{equation}

\subsection{Intermediate Jacobians and normal functions}
\label{subsec:normal-functions}

For each $y\in Y$, the intermediate Jacobian is the complex torus
\begin{equation}\label{eq:intermediate-jacobian}
    J^{2p-1}(X_y)
    :=\frac{\mathcal H_y}{F^p\mathcal H_y+\mathbb H_{\mathbb Z,y}}.
\end{equation}
These tori form a holomorphic family
$J^{2p-1}(\mathcal X/Y)\to Y$, whose sheaf of holomorphic sections
is locally described by
\[
    \frac{\mathcal H}{F^p\mathcal H+\mathbb H_{\mathbb Z}}.
\]

\begin{definition}
A \emph{normal function} is a holomorphic section
\[
    \nu\colon Y\longrightarrow J^{2p-1}(\mathcal X/Y)
\]
satisfying the horizontality condition, that is for every sufficiently small
open set $U\subset Y$ and every holomorphic lift
$\widetilde\nu\in H^0(U,\mathcal H)$ of $\nu|_U$,
\begin{equation}\label{eq:normal-function-transversality}
    \nabla\widetilde\nu
    \in H^0\bigl(U,F^{p-1}\mathcal H\otimes\Omega_Y^1\bigr).
\end{equation}
\end{definition}
This condition is independent of the lift (see \cite[Chapter~7]{VoisinII})

The normal functions considered in this paper arise from families
of algebraic cycles. Recall that, for a smooth projective variety
$X$, it is defined the Abel--Jacobi map
\[
    AJ_X:\mathrm{CH}^p(X)_{\mathrm{hom}}
    \longrightarrow J^{2p-1}(X),
\]
where $\mathrm{CH}^p(X)_{\mathrm{hom}}$ denotes the group of homologically trivial cycle classes of codimension $p$.
Let $\mathcal Z$ be a relative algebraic cycle of codimension
$p$ on $\mathcal X$, whose restriction $Z_y$ to each fiber is
homologically trivial. The fiberwise Abel--Jacobi classes define
a section
\[
    \nu_{\mathcal Z}:Y\longrightarrow
    J^{2p-1}(\mathcal X/Y),
    \qquad
    \nu_{\mathcal Z}(y):=AJ_{X_y}([Z_y]).
\]
This section is holomorphic and satisfies
\eqref{eq:normal-function-transversality}, hence defines a
normal function; see \cite[\S 7.2]{VoisinII}.

\subsection{Rank and canonical derivative}
\label{subsec:rank-canonical}

We use the description of the rank and canonical derivative of a normal
function developed by Hain \cite[\S\S 2.5, 3.1--3.2]{Hain}.

Recall that the Hodge structure on
$H^{2p-1}(X_y,\mathbb C)$ has weight $2p-1$.
Since this weight is odd, the Hodge decomposition gives
\[
    \mathcal H_y
    =
    F^p\mathcal H_y
    \oplus
    \overline{F^p\mathcal H_y}.
\]
The natural projection $H^{2p-1}(X_y,\mathbb R)
    \rightarrow
    \mathcal H_y/F^p\mathcal H_y$
is an isomorphism of real vector spaces. This yields an isomorphism of real tori
\begin{equation}\label{eq:real-jacobian-isomorphism}
    \frac{H^{2p-1}(X_y,\mathbb R)}
    {\mathbb H_{\mathbb Z,y}}
    \xrightarrow{\ \sim\ }J^{2p-1}(X_y).
\end{equation}

Thus the underlying real family of intermediate Jacobians is the
locally constant family
\[
    J^{2p-1}_{\mathbb R}(\mathcal X/Y)
    :=
    \frac{R^{2p-1}\pi_*\mathbb R}{\mathbb H_{\mathbb Z}}.
\]
Let $U\subset Y$ be a simply connected open set and fix
$y_0\in U$. The Gauss--Manin connection identifies the
cohomology groups of the fibers over $U$, giving a
trivialization
\[
    J^{2p-1}_{\mathbb R}(\mathcal X/Y)|_U
    \simeq
    \frac{H^{2p-1}(X_{y_0},\mathbb R)}
    {\mathbb H_{\mathbb Z,y_0}}\times U.
\]

For each fixed point $a$ of the real torus, the subsets $\mathcal F_a = \{a\}\times U$ define, through the identification
\eqref{eq:real-jacobian-isomorphism}, a foliation
$\mathcal F$ of $J^{2p-1}(\mathcal X/Y)$, called the
\emph{locally constant foliation}.
As observed by Hain \cite[\S 2.5]{Hain}, its leaves are
complex submanifolds. Indeed, locally constant cohomology
classes give flat sections of $\mathcal H$, whose images
in the intermediate Jacobian are holomorphic.
Thus $\mathcal F$ is a complex foliation, although it
need not be a holomorphic foliation.

Let $x\in J^{2p-1}(X_y)$, and let $\mathcal F_x$
be the leaf of $\mathcal F$ passing through $x$.
Since $\mathcal F_x$ is locally a section of
$J^{2p-1}(\mathcal X/Y)\to Y$, the differential of the
projection restricts to an isomorphism $T_x\mathcal F_x\simeq T_yY.$
On the other hand, the tangent space to the fiber at $x$
is naturally identified with $\mathcal H_y/F^p\mathcal H_y$.
We therefore obtain a direct sum decomposition
\[
    T_xJ^{2p-1}(\mathcal X/Y)
    =
    \frac{\mathcal H_y}{F^p\mathcal H_y}
    \oplus T_x\mathcal F_x.
\]
Equivalently, $T_x\mathcal F_x$ provides a splitting of the
tangent exact sequence
\[
    0\longrightarrow\mathcal H_y/F^p\mathcal H_y
    \longrightarrow T_xJ^{2p-1}(\mathcal X/Y)
    \longrightarrow T_yY\longrightarrow0.
\]
We denote by
\[
    \phi_x:T_xJ^{2p-1}(\mathcal X/Y)
    \longrightarrow\mathcal H_y/F^p\mathcal H_y
\]
the projection onto the vertical tangent space along
$T_x\mathcal F_x$. Since the leaves are complex submanifolds,
this projection is complex-linear.
Following Hain \cite[\S 3.1]{Hain}, we have the following definition.
\begin{definition}\label{def:rank-normal-function}
The \emph{rank} of a normal
function $\nu:Y\longrightarrow J^{2p-1}(\mathcal X/Y)$ at a point $y\in Y$ is defined by
\[
    \operatorname{rk}_y(\nu)
    :=
    \operatorname{rk}_{\mathbb C}
    \bigl(\phi_{\nu(y)}\circ d_y\nu\bigr).
\]
The \emph{rank of $\nu$} is
\[
    \operatorname{rk}(\nu)
    :=
    \max_{y\in Y}\operatorname{rk}_y(\nu).
\]
We say that $\nu$ has \emph{maximal rank} if $\operatorname{rk}(\nu)=\dim Y.$
\end{definition}

Now we describe the construction of the canonical derivative. Let $\nu_{\mathbb R}$ be the real-analytic section
of $J^{2p-1}_{\mathbb R}(\mathcal X/Y)$ corresponding to $\nu$.
Locally choose a real-analytic lift $\widetilde\nu_{\mathbb R}$
to the real cohomology bundle. Differentiating with its flat
connection gives a well-defined real-linear map
\begin{equation}\label{eq:real-derivative-normal-function}
    (\nabla\nu_{\mathbb R})_y:
    T_y^{\mathbb R}Y
    \longrightarrow H^{2p-1}(X_y,\mathbb R),
\end{equation}
where $T_y^{\mathbb R}Y$ is the tangent space of the underlying real
manifold. This does not depend on the lift, since two lifts differ
by a locally constant integral section.

Complexify \eqref{eq:real-derivative-normal-function} and restrict
to the holomorphic tangent space $T_yY=T_y^{1,0}Y$. This gives
\[
    (\nabla\nu_{\mathbb R})^{1,0}_y:
    T_yY\longrightarrow\mathcal H_y.
\]
We claim that the image of
$(\nabla\nu_{\mathbb R})^{1,0}_y$
is contained in $F^{p-1}\mathcal H_y$.
To see this, let $\widetilde\nu$ be a local holomorphic
lift of $\nu$ to $\mathcal H$.
Since $\widetilde\nu$ and $\widetilde\nu_{\mathbb R}$
represent the same section of the intermediate Jacobian,
we may arrange, by adding a locally constant integral
section to one of the lifts, that
\[
    \eta:=\widetilde\nu-\widetilde\nu_{\mathbb R}
\]
takes values in $F^p\mathcal H$.
Although $\eta$ is only real-analytic, Griffiths
transversality implies that its $(1,0)$--part of the derivative takes values in
$F^{p-1}\mathcal H\otimes\Omega_Y^1$.
Indeed, this follows by expressing $\eta$ locally
in a holomorphic frame of $F^p\mathcal H$.
On the other hand, the defining condition of a normal function gives $\nabla\widetilde\nu
    \in F^{p-1}\mathcal H\otimes\Omega_Y^1.$
Therefore,
\[
    (\nabla\widetilde\nu_{\mathbb R})^{1,0}
    =
    \nabla\widetilde\nu-\nabla^{1,0}\eta
    \in
    F^{p-1}\mathcal H\otimes\Omega_Y^1,
\]
as claimed.

The \emph{canonical derivative} of $\nu$ at $y$, in Hain's terminology
\cite[\S 3.2]{Hain}, is therefore the complex-linear map
\begin{equation}\label{eq:canonical-derivative-definition}
    \overline\nabla_y\nu:
    T_yY\longrightarrow
    F^{p-1}\mathcal H_y/F^p\mathcal H_y,
    \qquad
    \xi\longmapsto
    \bigl[(\nabla\nu_{\mathbb R})_y^{1,0}(\xi)\bigr].
\end{equation}
It depends real-analytically on $y$, but is not in general a
holomorphic section of
\[
    (F^{p-1}\mathcal H/F^p\mathcal H)
    \otimes\Omega_Y^1.
\]

The preceding constructions satisfy the following rank identity,
which is the odd-weight cohomological form of
\cite[Proposition~3.3]{Hain}.

\begin{proposition}\label{prop:rank-canonical-derivative}
For every $y\in Y$,
\[
    \operatorname{rk}_y(\nu)
    =
    \operatorname{rk}_{\mathbb C}(\overline\nabla_y\nu)
    =
    \frac12\operatorname{rk}_{\mathbb R}
       \bigl((\nabla\nu_{\mathbb R})_y\bigr).
\]
\end{proposition}

\begin{proof}
In a local trivialization induced by the Gauss--Manin
connection, $\nu_{\mathbb R}$ is a map to a fixed real
torus, and $(\nabla\nu_{\mathbb R})_y$ is its differential.
Under the identification
$H^{2p-1}(X_y,\mathbb R)\simeq\mathcal H_y/F^p\mathcal H_y$,
this agrees with the vertical component
$\phi_{\nu(y)}\circ d_y\nu$.
The $(1,0)$ component of this map, expressed in the quotient
$\mathcal H_y/F^p\mathcal H_y$, is precisely
$\overline\nabla_y\nu$. Since the real rank of a complex-linear
map is twice its complex rank, the result follows.
\end{proof}

\begin{remark}\label{remark dense open}
    As observed by Hain \cite[Proposition~3.6]{Hain}, the rank is lower
semicontinuous. Indeed, the canonical derivative is real analytic, so
for each $r\geq0$ the locus
\[
    Y_{\leq r}:=\{y\in Y: \operatorname{rk}_y(\nu)\leq r\}
\]
is a real-analytic subset. If $Y$ is connected of dimension $d$ and
$\operatorname{rk}_{y_0}(\nu)=d$ at some $y_0\in Y$, then
\[
    \{y\in Y:\operatorname{rk}_y(\nu)=d\}
    =Y\setminus Y_{\leq d-1}
\]
is open and dense in $Y$. Thus maximal rank at one point implies
maximal rank on a dense open subset in the analytic topology.
\end{remark}

\subsection{The infinitesimal invariant}
\label{subsec:infinitesimal-invariant}

The infinitesimal invariant of a normal function was introduced by Griffiths \cite{Griffiths1983} and subsequently studied by Green \cite{Green1989} and Voisin \cite{Voisin1988}.
We recall its construction in the form used by
Collino--Naranjo--Pirola \cite[\S 5.1]{CNP}.

By Griffiths transversality, the Gauss--Manin connection
induces an $\mathcal O_Y$-linear morphism
\begin{equation}\label{eq:IVHS-p}
    \nabla_p^{p-1}:
    \frac{F^p\mathcal H}{F^{p+1}\mathcal H}
    \longrightarrow
    \frac{F^{p-1}\mathcal H}{F^p\mathcal H}
    \otimes\Omega_Y^1.
\end{equation}

Let $\nu$ be a normal function and choose a local holomorphic
lift $\widetilde\nu$ to $\mathcal H$. By horizontality \eqref{eq:normal-function-transversality},
$\nabla\widetilde\nu$ takes values in
$F^{p-1}\mathcal H\otimes\Omega_Y^1$.
Changing the lifting $\widetilde\nu$ changes
$\nabla\widetilde\nu$ modulo $F^p\mathcal H$ by an element
of $\operatorname{Im}(\nabla_p^{p-1})$.
Thus its class modulo this image is independent of the
choice of lift; see \cite[\S 7.1.3]{VoisinII}.

\begin{definition}\label{def:infinitesimal-invariant}
Following \cite[\S 5.1]{CNP}, the \emph{infinitesimal
invariant} of $\nu$ is the section
\[
    \delta(\nu)\in
    H^0\bigl(Y,\operatorname{Coker}(\nabla_p^{p-1})\bigr)
\]
locally represented by $\nabla\widetilde\nu$ modulo
$F^p\mathcal H$. Its value at $y\in Y$ will be denoted
by $\delta_y(\nu)$.
\end{definition}

We now relate the infinitesimal invariant to the
canonical derivative. Choose local
holomorphic and real-analytic lifts
$\widetilde\nu$ and $\widetilde\nu_{\mathbb R}$ of $\nu$
as in Section~\ref{subsec:rank-canonical}. After adjusting
by a flat integral section, their difference
\[
    \eta:=\widetilde\nu-\widetilde\nu_{\mathbb R}
\]
is a real-analytic section of $F^p\mathcal H$.
Taking $(1,0)$ derivatives and reducing modulo
$F^p\mathcal H_y$ gives
\[
    [(\nabla\widetilde\nu)_y]
    =
    \overline\nabla_y\nu
    +
    (\nabla_p^{p-1})_y([\eta_y]).
\]
Here $[\eta_y]$ denotes the class of $\eta_y$ in
$F^p\mathcal H_y/F^{p+1}\mathcal H_y$.
The formula remains valid even though $\eta$ is not
holomorphic: in a local holomorphic frame of
$F^p\mathcal H$, the derivatives of its coefficients
contribute only terms in $F^p\mathcal H$.
Consequently we have
\begin{equation}\label{eq:delta-canonical}
    \delta_y(\nu)
    =
    [\overline\nabla_y\nu]
    \in
    \operatorname{Coker}\bigl((\nabla_p^{p-1})_y\bigr).
\end{equation}
In particular we have the relation
\begin{equation}\label{eq:delta-nonzero-canonical}
    \delta_y(\nu)\neq0
    \quad\Longrightarrow\quad
    \overline\nabla_y\nu\neq0.
\end{equation}

% -------------------------------------------------
\subsection{Symmetries}
\label{subsec:symmetries}
% -------------------------------------------------

We conclude with the equivariance principle for the canonical
derivative, following \cite[Lemma~3.5]{Hain}.

Let $G$ be a finite group acting holomorphically on $Y$,
and suppose that its action lifts to $\mathcal X$
compatibly with $\pi$.
The induced action on $\mathbb H_{\mathbb Z}$ and
$\mathcal H$ preserves the Gauss--Manin connection
and the Hodge filtration.
Both Jacobian bundles therefore inherit a $G$-action,
and their natural identification is $G$-equivariant.

A normal function $\nu$ is said to be \emph{$G$-invariant} if it is
equivariant as a section, that is,
\[
    \nu(g\cdot y)=g\cdot\nu(y),
    \qquad g\in G,\quad y\in Y.
\]

\begin{proposition}
\label{prop:symmetry-canonical-derivative}
Let $\nu$ be a $G$-invariant normal function and let $y\in Y$
be fixed by $G$. Then
\[
    \overline\nabla_y\nu
    \in
    \operatorname{Hom}_G\left(
        T_yY,\,
        F^{p-1}\mathcal H_y/F^p\mathcal H_y
    \right).
\]
\end{proposition}

\begin{proof}
The induced real-analytic section $\nu_{\mathbb R}$ is
$G$-equivariant. Since the action on the real cohomology bundle
preserves its flat connection, differentiating at the fixed point
$y$ gives
\[
    (\nabla\nu_{\mathbb R})_y\circ d_yg
    =
    g\circ(\nabla\nu_{\mathbb R})_y,
    \qquad g\in G.
\]
Here $g$ on the right denotes the induced action on
$H^{2p-1}(X_y,\mathbb R)$; the identity is independent of the
choice of a local real lift.
Since $G$ acts holomorphically on $Y$, passage to the $(1,0)$
component preserves equivariance. The Hodge filtration is also
$G$-stable, so projection to
$F^{p-1}\mathcal H_y/F^p\mathcal H_y$ is equivariant.
By \eqref{eq:canonical-derivative-definition}, this projection
of $(\nabla\nu_{\mathbb R})_y^{1,0}$ is precisely
$\overline\nabla_y\nu$.
\end{proof}

%------------------------------------------------
\section{The Fano normal function}
\label{sec:fano-normal-function}
%------------------------------------------------

We recall the construction of the Fano cycle and identify the
primitive component of its Abel--Jacobi class that will be used below.

\subsection{The Fano cycle}
\label{subsec:fano-cycle}

Let $V=\{E=0\}\subset\mathbb P^4$ be a smooth cubic threefold and
$F=F(V)$ its Fano surface of lines. The intermediate Jacobian
\[
    J(V):=
    \frac{H^3(V,\mathbb C)}
         {F^2H^3(V,\mathbb C)+H^3(V,\mathbb Z)}
\]
is a principally polarized abelian variety of dimension $5$; we denote
its polarization class by $\theta\in H^{1,1}(J(V))\cap H^2(J(V),\mathbb Z)$.

The Abel--Jacobi map for lines induces the classical identification
$\Alb(F)\simeq J(V)$ \cite{ClemensGriffiths,CNP}. Choosing a base
line $r\in F$ gives an Albanese embedding
\[
    a_r:F\hookrightarrow J(V),\qquad
    q\longmapsto\AJ_V(\ell_q-\ell_r),
\]
where $\ell_q\subset V$ is the line represented by $q\in F$.
Set $F_r:=a_r(F)$ and $F_r^-:=[-1](F_r)$, where $[-1]$ denotes
involution on $J(V)$.

Since $[-1]$ acts on $H^m(J(V),\mathbb Z)$ as multiplication by
$(-1)^m$, the surfaces $F_r$ and $F_r^-$ have the same class in
$H^6(J(V),\mathbb Z)$. Thus
\begin{equation}\label{eq:Fano-cycle}
    Z_r:=F_r-F_r^-
\end{equation}
is a homologically trivial algebraic cycle of codimension $3$ on
$J(V)$, called the \emph{Fano cycle} \cite[\S 5.2]{CNP}.
Its Abel--Jacobi class is
\begin{equation}\label{eq:nu-r}
    \nu_r(V):=\AJ_{J(V)}(Z_r)\in J^5(J(V)),
\end{equation}
where, according to \eqref{eq:intermediate-jacobian},
\[
    J^5(J(V))=
    \frac{H^5(J(V),\mathbb C)}
         {F^3H^5(J(V),\mathbb C)+H^5(J(V),\mathbb Z)}.
\]

Changing the base line translates $F_r$ inside $J(V)$, so the class
$\nu_r(V)$ may depend on $r$.

%------------------------------------------------
\subsection{The primitive fifth cohomology}
\label{subsec:primitive-fifth}
%------------------------------------------------

We follow the Lefschetz description of the primitive Fano invariant
in \cite[\S 5.3]{CNP}. Put $A=J(V)$ and $L:=\theta\wedge-$. The
Lefschetz decomposition gives an isomorphism of rational Hodge structures
\begin{equation}\label{eq:lefschetz-H5}
    H^5(A,\mathbb Q)
    =\theta^2H^1(A,\mathbb Q)
     \oplus\theta P^3(A,\mathbb Q)
     \oplus P^5(A,\mathbb Q),
\end{equation}
where
\[
    P^3(A,\mathbb Q):=\ker(L^3:H^3(A,\mathbb Q)\to H^9(A,\mathbb Q)),
    \qquad
    P^5(A,\mathbb Q):=\ker(L:H^5(A,\mathbb Q)\to H^7(A,\mathbb Q)).
\]
Since $\theta$ has type $(1,1)$, the primitive fifth cohomology is a
rational Hodge substructure, with
\[
    P^5(A,\mathbb C)=\bigoplus_{a+b=5}P^{a,b}(A).
\]

Let $\operatorname{pr}_5:H^5(A,\mathbb Q)\to P^5(A,\mathbb Q)$
be the projection, and set
$\Lambda_5:=\operatorname{pr}_5H^5(A,\mathbb Z)$. We define
\begin{equation}\label{eq:primitive-J5}
    J^5_{\mathrm{prim}}(A):=
    \frac{P^5(A,\mathbb C)}
    {F^3P^5(A,\mathbb C)+\Lambda_5}.
\end{equation}
Since $\operatorname{pr}_5$ preserves the Hodge filtration and
maps $H^5(A,\mathbb Z)$ onto $\Lambda_5$, it induces a
well-defined homomorphism
\[
    \operatorname{pr}_5:J^5(A)\longrightarrow
    J^5_{\mathrm{prim}}(A).
\]

We recall the following result of Collino--Naranjo--Pirola.

\begin{proposition}[{\cite[Proposition~5.4]{CNP}}]
\label{prop:CNP-primitive}
The component of $\nu_r(V)$ associated with $\theta P^3(A)$
vanishes, and its primitive fifth component is independent
of the base line $r\in F(V)$.
\end{proposition}

Thus the primitive Fano invariant is well defined by
\begin{equation}\label{eq:primitive-AJ-Fano}
    \nu_5(V):=\operatorname{pr}_5(\nu_r(V))
    \in J^5_{\mathrm{prim}}(J(V)).
\end{equation}

%------------------------------------------------
\subsection{The Fano normal function in families}
\label{subsec:fano-normal-function}
%------------------------------------------------

Let $\pi:\mathcal V\to Y$ be a smooth family of cubic threefolds,
and let $\mathcal J\to Y$ be the associated family of intermediate
Jacobians. Locally on $Y$, a section of the relative Fano surface
gives a family of homologically trivial Fano cycles in
$\mathcal J$. Their Abel--Jacobi classes define a normal function
\[
    \nu_r:Y\longrightarrow J^5(\mathcal J/Y).
\]
The cohomology class of the principal polarization is preserved by the Gauss--Manin connection. Consequently, the primitive cohomology groups
$P^5(J(V_y),\mathbb Q)$ form a variation of Hodge
structure over $Y$.
By Proposition~\ref{prop:CNP-primitive}, the primitive component of
$\nu_r$ is independent of the choice of the local base line. These
components therefore glue to a normal function
\[
    \nu_F:Y\longrightarrow J^5_{\mathrm{prim}}(\mathcal J/Y),
\]
called the \emph{Fano normal function}.
Its underlying variation of Hodge structure is the primitive fifth
cohomology. Thus the canonical derivative introduced in
\eqref{eq:canonical-derivative-definition} has the form
\begin{equation}\label{eq:canonical-Fano-general}
    \overline\nabla_y\nu_F:
    T_yY\longrightarrow
    \frac{F^2P^5(J(V_y),\mathbb C)}
         {F^3P^5(J(V_y),\mathbb C)}
    \simeq P^{2,3}(J(V_y)).
\end{equation}

%------------------------------------------------
\subsection{The Jacobian ring and the relevant Hodge spaces}
\label{subsec:jacobian-ring}
%------------------------------------------------

Let $V=\{E=0\}\subset\mathbb P^4$ be a smooth cubic threefold.
Its Jacobian ideal and Jacobian ring are
\[
    J_E:=\left(\frac{\partial E}{\partial x_0},\ldots,
                 \frac{\partial E}{\partial x_4}\right),
    \qquad
    R^\bullet(V):=
    \frac{\mathbb C[x_0,\ldots,x_4]}{J_E}
    =\bigoplus_{d\geq0}R^d(V).
\]
The nonzero graded pieces have dimensions
\[
    \bigl(\dim R^d(V)\bigr)_{d=0}^{5}
    =(1,5,10,10,5,1),
\]
and multiplication induces perfect pairings
\[
    R^d(V)\otimes R^{5-d}(V)\longrightarrow R^5(V)\simeq\mathbb C.
\]

By Griffiths' residue theorem
\cite[\S 6.1.3]{VoisinII}, the Hodge components of
a smooth cubic threefold $V$ can be identified with
graded pieces of its Jacobian ring. After choosing
a volume form, we have
\begin{equation}\label{eq:jacobian-hodge}
    H^{2,1}(V)\simeq R^1(V),
    \qquad
    H^{1,2}(V)\simeq R^4(V).
\end{equation}

Let $(Y,0)$ be a local deformation space of $V$, chosen
so that the Kodaira--Spencer map is an isomorphism at $0$.
We then have
\begin{equation}\label{eq:tangent-R3}
    T_0Y\simeq H^1(V,T_V)\simeq R^3(V).
\end{equation}
Under these identifications, the infinitesimal variation
of Hodge structure is given by multiplication in the
Jacobian ring:
\begin{equation}\label{eq:IVHS-jacobian-ring}
    R^3(V)\otimes R^1(V)\longrightarrow R^4(V).
\end{equation}
See \cite[\S 3]{CNP}.

The incidence correspondence for lines on $V$, together with
$\operatorname{Alb}(F(V))\simeq J(V)$, identifies
\[
    H^{1,0}(J(V))\simeq H^{2,1}(V),\qquad
    H^{0,1}(J(V))\simeq H^{1,2}(V)
\]
(see \cite{ClemensGriffiths}). Set
\begin{equation}\label{eq:MN-definition}
    M:=H^{1,0}(J(V))\simeq R^1(V),\qquad
    N:=H^{0,1}(J(V))\simeq R^4(V).
\end{equation}
Since the cohomology algebra of an abelian variety is the exterior
algebra of its first cohomology, we have
\[
    H^{2,3}(J(V))\simeq\Lambda^3N\otimes\Lambda^2M,\qquad
    H^{1,2}(J(V))\simeq\Lambda^2N\otimes M.
\]
In bidegree $(2,3)$, the Lefschetz decomposition gives
\[
    H^{2,3}(J(V))
    =P^{2,3}(J(V))\oplus
      \theta\wedge H^{1,2}(J(V)).
\]
Consequently,
\begin{equation}\label{eq:primitive-target}
    P^{2,3}(J(V))\simeq
    \frac{\Lambda^3N\otimes\Lambda^2M}
         {\theta\wedge(\Lambda^2N\otimes M)}.
\end{equation}

Combining \eqref{eq:canonical-Fano-general},
\eqref{eq:tangent-R3} and \eqref{eq:primitive-target},
we can express the canonical derivative of the Fano
normal function at $V$ as
\begin{equation}\label{eq:Fano-derivative-algebraic}
    \overline\nabla_V\nu_F:
    R^3(V)\longrightarrow
    \frac{\Lambda^3N\otimes\Lambda^2M}
         {\theta\wedge(\Lambda^2N\otimes M)}.
\end{equation}

%------------------------------------------------
\subsection{Symmetries of the Fano normal function}
\label{subsec:fano-symmetries}
%------------------------------------------------

We verify that the Fano normal function satisfies the equivariance
hypothesis of Proposition~\ref{prop:symmetry-canonical-derivative}.
Let $\pi:\mathcal V\to Y$ be a smooth family of cubic threefolds
with a finite group $G$ acting on the family by projective
isomorphisms, covering a holomorphic action on $Y$. Thus each
$g\in G$ induces an isomorphism $V_t\to V_{g\cdot t}$.
By functoriality, $G$ acts on the associated family $\mathcal J$ of
intermediate Jacobians, preserving their principal polarizations.
In particular, it preserves the primitive fifth cohomology and the
corresponding family $J^5_{\mathrm{prim}}(\mathcal J/Y)$.

\begin{proposition}\label{prop:Fano-G-invariant}
The Fano normal function is $G$-equivariant:
\[
    \nu_F(g\cdot t)=g\cdot\nu_F(t)
    \qquad(g\in G,\ t\in Y).
\]
\end{proposition}

\begin{proof}
Fix $t\in Y$ and a base line $r\in F(V_t)$. The isomorphism
$g:V_t\to V_{g\cdot t}$ induces isomorphisms $g_F$ of Fano
surfaces and $g_J$ of their intermediate Jacobians. Functoriality
of the Abel--Jacobi map for lines gives
\[
    g_J\circ a_r=a_{g_F(r)}\circ g_F.
\]
As $g_J$ commutes with the involution, it sends the Fano cycle
$Z_r=F_r-F_r^-$ to $Z_{g_F(r)}$. Functoriality of Abel--Jacobi
classes, followed by the polarization-preserving primitive
projection, therefore yields
\[
    g\cdot\nu_{5,r}(V_t)
    =\nu_{5,g_F(r)}(V_{g\cdot t}).
\]
The primitive component is independent of the base line by Proposition \ref{prop:CNP-primitive}. Thus the right-hand side is $\nu_F(g\cdot t)$.
\end{proof}

If $t_0\in Y$ is fixed by $G$, the Proposition \ref{prop:Fano-G-invariant} and
Proposition~\ref{prop:symmetry-canonical-derivative} give
\begin{equation}\label{eq:Fano-canonical-equivariant-general}
    \overline\nabla_{t_0}\nu_F\in
    \operatorname{Hom}_G\bigl(
        T_{t_0}Y,\,P^{2,3}(J(V_{t_0}))
    \bigr).
\end{equation}
This is the equivariance property used in the computation at the
Fermat cubic.

%------------------------------------------------
\section{The Fermat cubic threefold}
\label{sec:fermat}
%------------------------------------------------

We now specialize to the Fermat cubic threefold
\[
    V_F=\{E_F=x_0^3+\cdots+x_4^3=0\}
    \subset\mathbb P^4.
\]
Its Jacobian ring is
\[
    R_F^\bullet
    =
    \frac{\mathbb C[x_0,\ldots,x_4]}
         {(x_0^2,\ldots,x_4^2)},
\]
so that $R_F^d$ has a basis of square-free monomials of degree $d$.

The automorphism group of the Fermat cubic is
\[
    G=\operatorname{Aut}(V_F)=A\rtimes S_5,
    \qquad
    A:=\left\{
    (\lambda_0,\ldots,\lambda_4)\in\mu_3^5:
    \prod_{i=0}^4\lambda_i=1
    \right\}
    \simeq(\mu_3)^4;
\]
see \cite[Example~3.1(1)]{WeiYu2020}.
The group is generated by diagonal transformations
and permutations of the coordinates. We write
$g=(\lambda_0,\ldots,\lambda_4;\sigma)$ and fix the convention
\begin{equation}\label{eq:fermat-polynomial-action}
    g\cdot x_i=\lambda_i x_{\sigma(i)},
    \qquad \sigma\in S_5.
\end{equation}
Here the action on equations is induced by inverse pullback; we use
the corresponding action on cohomology, so that these are left
$G$-representations.

Let $(Y,0)$ be a local deformation space of $V_F$,
chosen to be $G$-invariant and such that the
Kodaira--Spencer map is an isomorphism at the origin.
The action of $G$ on $V_F$ extends to the family,
and we have $T_0Y\simeq R_F^3$. By Proposition~\ref{prop:Fano-G-invariant}, the Fano normal function is
$G$-equivariant; therefore
\begin{equation}\label{eq:Fano-derivative-Fermat-equivariant}
    \overline\nabla_0\nu_F\in
    \operatorname{Hom}_G\bigl(T_0Y,P^{2,3}(J(V_F))\bigr).
\end{equation}

%------------------------------------------------
\subsection{The action on the tangent and Hodge spaces}
\label{subsec:fermat-Hodge-action}
%------------------------------------------------

Let $\chi:G\to\{\pm1\}$ be the determinant character of the
chosen linear action on $\mathbb C^5$. Since
$\prod_i\lambda_i=1$,
\begin{equation}\label{eq:det-character}
    \chi(g)=\operatorname{sgn}(\sigma).
\end{equation}

\begin{lemma}\label{lem:fermat-tangent-Hodge}
With the action \eqref{eq:fermat-polynomial-action}, we have isomorphisms of $G$-representations
\[
    T_0Y\simeq R_F^3,\qquad
    M\simeq R_F^1\otimes\chi,\qquad
    N\simeq R_F^4\otimes\chi.
\]
\end{lemma}

\begin{proof}

Since $G$ preserves the Fermat polynomial $E_F$ and its
Jacobian ideal, its natural action on cubic polynomials
induces an action on $R_F^3$.
Under the Kodaira--Spencer identification
$T_0Y\simeq R_F^3$, this agrees with the action on
first-order deformations. In particular, this
identification is $G$-equivariant, and no character
twist occurs.

For the Hodge spaces, consider the homogeneous form
\[
    \Omega_{\mathbb P}
    :=\sum_{i=0}^4(-1)^i x_i\,
       dx_0\wedge\cdots\wedge\widehat{dx_i}
       \wedge\cdots\wedge dx_4.
\]
It is the contraction of $dx_0\wedge\cdots\wedge dx_4$ with the
Euler vector field; consequently, a linear coordinate transformation
acts on it by its determinant. Thus
\[
    g\cdot\Omega_{\mathbb P}
    =\chi(g)\Omega_{\mathbb P}.
\]

By \eqref{eq:jacobian-hodge}, the residue maps identify
$R_F^1$ and $R_F^4$ with $H^{2,1}(V_F)$ and
$H^{1,2}(V_F)$, respectively.
More precisely, these identifications are given by
\[
    P\longmapsto
    \operatorname{Res}_{V_F}
    \left(\frac{P\Omega_{\mathbb P}}{E_F^2}\right),
    \qquad
    Q\longmapsto
    \operatorname{Res}_{V_F}
    \left(\frac{Q\Omega_{\mathbb P}}{E_F^3}\right),
\]
where the residue classes are taken in the corresponding
graded pieces of the Hodge filtration
\cite[\S 6.1.3]{VoisinII}.
The action of $G$ preserves $E_F$, while the homogeneous
volume form satisfies
$g\cdot\Omega_{\mathbb P}=\chi(g)\Omega_{\mathbb P}$.
Consequently,
\[
    g\cdot\operatorname{Res}_{V_F}
    \left(\frac{P\Omega_{\mathbb P}}{E_F^q}\right)
    =
    \chi(g)\operatorname{Res}_{V_F}
    \left(\frac{(g\cdot P)\Omega_{\mathbb P}}{E_F^q}\right),
    \qquad q=2,3.
\]
Thus, using \eqref{eq:MN-definition}, the thesis follows.

\end{proof}

%------------------------------------------------
\subsection{Bases and weights}
\label{subsec:bases-weights}
%------------------------------------------------

We choose bases of the $G$-representations
$M\simeq R_F^1\otimes\chi$ and $N\simeq R_F^4\otimes\chi$.
Set
\[
    e_i:=x_i\in R_F^1,\qquad
    f_i:=\prod_{j\ne i}x_j\in R_F^4,
    \qquad i=0,\ldots,4.
\]
The vectors $e_i$ and $f_i$ form bases of $R_F^1$ and $R_F^4$,
respectively. For $g=(\lambda_0,\ldots,\lambda_4;\sigma)\in G$,
the polynomial action satisfies
\begin{equation}\label{eq:action-ei-fi-polynomial}
    g\cdot e_i=\lambda_i e_{\sigma(i)},\qquad
    g\cdot f_i=\lambda_i^{-1}f_{\sigma(i)},
\end{equation}
since $\prod_{j\ne i}\lambda_j=\lambda_i^{-1}$.
Using the same symbols for the corresponding vectors in $M$ and
$N$, Lemma~\ref{lem:fermat-tangent-Hodge} gives
\begin{equation}\label{eq:action-B-E}
    g\cdot e_i=\chi(g)\lambda_i e_{\sigma(i)},\qquad
    g\cdot f_i=\chi(g)\lambda_i^{-1}f_{\sigma(i)}.
\end{equation}

We next describe the weights of the diagonal subgroup $A$.
For $i=0,\ldots,4$, let $\varepsilon_i\in A^\vee$ be the
coordinate character
\[
    \varepsilon_i(\lambda_0,\ldots,\lambda_4)=\lambda_i.
\]
We write the character group additively. The relation
$\prod_i\lambda_i=1$ gives $\varepsilon_0+\cdots+\varepsilon_4=0.$
More precisely,
\begin{equation}\label{eq:character-group-A}
    A^\vee\simeq
    \frac{(\mathbb Z/3\mathbb Z)^5}
         {\langle(1,1,1,1,1)\rangle},
\end{equation}
where $\varepsilon_i$ is represented by the $i$-th standard
basis vector. Since $\chi|_A=1$, the determinant twist does not
change the $A$-weights:
\[
    \operatorname{wt}(e_i)=\varepsilon_i,\qquad
    \operatorname{wt}(f_i)=-\varepsilon_i.
\]
For any $S\subset\{0,\ldots,4\}$, we shall write $\varepsilon_S:=\sum_{i\in S}\varepsilon_i.$

%------------------------------------------------

\subsection{The Lefschetz operator and the primitive quotient}
\label{subsec:lefschetz-fermat}

Let $\theta\in H^{1,1}(J(V_F))$ be the class of the
principal polarization. Recall that this polarization
is induced by the intersection pairing on
$H^3(V_F,\mathbb Z)$ \cite{ClemensGriffiths}.
Since every automorphism of $V_F$ preserves this
pairing, the class $\theta$ is $G$-invariant.
We now determine its expression in terms of the
bases of $M$ and $N$ introduced above.

\begin{lemma}\label{lem:theta-invariant}
Under the identification
$H^{1,1}(J(V_F))\simeq N\otimes M$, we have
\begin{equation}\label{eq:theta-Fermat}
    \theta=c\sum_{i=0}^4 f_i\otimes e_i
\end{equation}
for some $c\in\mathbb C^*$.
\end{lemma}

\begin{proof}
Since $\theta\in N\otimes M$, we can write
\[
    \theta=\sum_{i,j=0}^4 c_{ij}f_i\otimes e_j.
\]
Recall that $f_i\otimes e_j$ has $A$-weight
$-\varepsilon_i+\varepsilon_j$, which is trivial
if and only if $i=j$. Since $\theta$ is $A$-invariant,
it follows that
\[
    \theta=\sum_{i=0}^4 c_i f_i\otimes e_i.
\]
Moreover, by \eqref{eq:action-B-E}, the action of
$S_5$ permutes the vectors $f_i\otimes e_i$.
As $\theta$ is also $S_5$-invariant, all coefficients
$c_i$ must be equal. Hence
\[
    \theta=c\sum_{i=0}^4 f_i\otimes e_i
\]
for some $c\in\mathbb C$. Finally, $c\neq 0$,
since $\theta$ is the class of the principal polarization.
\end{proof}
Since multiplication by a nonzero scalar does not
change the image of the Lefschetz operator, we may
write
\[
    \theta = \sum_{i=0}^4 f_i\otimes e_i
\]
Wedge product with $\theta$ gives the Lefschetz map
\begin{equation}\label{eq:Lefschetz-Fermat}
    L_\theta:\Lambda^2N\otimes M
    \longrightarrow\Lambda^3N\otimes\Lambda^2M,
    \qquad
    u\otimes b\longmapsto
    \sum_{i=0}^4(f_i\wedge u)\otimes(e_i\wedge b).
\end{equation}

\begin{lemma}\label{lem:Lefschetz-equivariant}
The map $L_\theta$ is $G$-equivariant.
\end{lemma}

\begin{proof}
For $g\in G$ and $v\in\Lambda^2N\otimes M$, the invariance of
$\theta$ gives
\[
    g\cdot L_\theta(v)
    =g\cdot(\theta\wedge v)
    =\theta\wedge(g\cdot v)
    =L_\theta(g\cdot v).
\]
\end{proof}
Thus $\operatorname{Im}(L_\theta)$ is $G$-stable, and the primitive
quotient
\begin{equation}\label{eq:Q-Fermat}
    Q:=\frac{\Lambda^3N\otimes\Lambda^2M}
    {\theta\wedge(\Lambda^2N\otimes M)}
\end{equation}
is a $G$-representation. By \eqref{eq:primitive-target},
$Q\simeq P^{2,3}(J(V_F))$ as $G$-representations. Hence
\begin{equation}\label{eq:canonical-HomG-Fermat}
    \overline\nabla_0\nu_F\in\operatorname{Hom}_G(R_F^3,Q).
\end{equation}

%------------------------------------------------
\subsection{Weights in the domain and in the primitive quotient}
\label{subsec:weights-domain-codomain}
%------------------------------------------------

We compare the $A$-weights of $R_F^3$ with those of the primitive
quotient $Q$. For $S\subset\{0,\ldots,4\}$ with $|S|=3$, put
\[
    v_S:=\prod_{i\in S}x_i\in R_F^3.
\]
The ten square-free monomials $v_S$ form a basis of $R_F^3$, and
\[
    \operatorname{wt}(v_S)=\varepsilon_S,\qquad
    R_F^3=\bigoplus_{|S|=3}\mathbb C v_S.
\]

For any subset $I\subset\{0,\ldots,4\}$, let
$\mathbf 1_I\in(\mathbb Z/3\mathbb Z)^5$ be its
\emph{indicator vector}, defined by
\[
    (\mathbf 1_I)_i=
    \begin{cases}
        1 & \text{if } i\in I,\\
        0 & \text{if } i\notin I.
    \end{cases}
\]
By \eqref{eq:character-group-A}, the character
$\varepsilon_I$ is represented by $\mathbf 1_I$,
and two vectors represent the same character of $A$
if and only if their difference is a multiple of
$\mathbf 1=(1,1,1,1,1)$.
In particular, the ten characters $\varepsilon_S$,
with $|S|=3$, are pairwise distinct.
Indeed, if $\varepsilon_S=\varepsilon_T$, then $\mathbf 1_S-\mathbf 1_T=c\mathbf 1$
for some $c\in\mathbb Z/3\mathbb Z$.
Summing the coordinates and using $|S|=|T|=3$,
we obtain $5c=0$ in $\mathbb Z/3\mathbb Z$.
Thus $c=0$, and hence $S=T$.

For any subsets $I,J\subset\{0,\ldots,4\}$, write
\[
    f_I:=\bigwedge_{i\in I}f_i,\qquad
    e_J:=\bigwedge_{j\in J}e_j,
\]
with increasing order of indices. If $|I|=3$ and $|J|=2$,
the tensors $f_I\otimes e_J$ form a basis of
$\Lambda^3N\otimes\Lambda^2M$, with
\[
    \operatorname{wt}(f_I\otimes e_J)
    =-\varepsilon_I+\varepsilon_J.
\]
For every $|S|=3$, define
\begin{equation}\label{eq:qS-definition}
    q_S:=f_S\otimes e_{S^c}
    \in\Lambda^3N\otimes\Lambda^2M.
\end{equation}

\begin{lemma}\label{lem:numerator-weight-space}
For every $|S|=3$, the $\varepsilon_S$-weight space of the
numerator of $Q$ is
\[
    (\Lambda^3N\otimes\Lambda^2M)_{\varepsilon_S}
    =\mathbb C q_S.
\]
\end{lemma}

\begin{proof}
Since $\varepsilon_{S^c}=-\varepsilon_S$ and the character
group of $A$ has exponent $3$, we have
\[
    \operatorname{wt}(q_S)
    =-\varepsilon_S+\varepsilon_{S^c}
    =-2\varepsilon_S=\varepsilon_S.
\]

Conversely, suppose that a basis vector $f_I\otimes e_J$,
with $|I|=3$ and $|J|=2$, has weight $\varepsilon_S$. Then
\[
    -\mathbf 1_I+\mathbf 1_J-\mathbf 1_S
    =c\mathbf 1
\]
for some $c\in\mathbb Z/3\mathbb Z$. Summing coordinates yields
\[
    -3+2-3\equiv5c\pmod3,
\]
so $c=1$. If $i\in S$, the $i$-th coordinate now gives
\[
    -\mathbf 1_I(i)+\mathbf 1_J(i)\equiv-1\pmod3.
\]
This forces
$\mathbf 1_I(i)=1$ and $\mathbf 1_J(i)=0$. Similarly, for
$i\notin S$ one obtains
$\mathbf 1_I(i)=0$ and $\mathbf 1_J(i)=1$.
Therefore $I=S$ and $J=S^c$.
\end{proof}

\begin{lemma}\label{lem:Lefschetz-no-epsilonS}
For every $|S|=3$,
\[
    (\Lambda^2N\otimes M)_{\varepsilon_S}=0,
    \qquad
    Q_{\varepsilon_S}=\mathbb C[q_S],
\]
where $[q_S]$ is the image of $q_S$ in $Q$.
\end{lemma}

\begin{proof}
A basis vector of $\Lambda^2N\otimes M$ is
$f_I\otimes e_k$, with $|I|=2$, and has weight
$-\varepsilon_I+\varepsilon_k$. If this weight were
$\varepsilon_S$, we would have
\[
    -\mathbf 1_I+\mathbf 1_{\{k\}}-\mathbf 1_S
    =c\mathbf 1.
\]
Summing coordinates gives
$-2+1-3\equiv5c\pmod3$, hence $c=1$.
For every $i\in S$, the $i$-th coordinate would then imply
\[
    -\mathbf 1_I(i)+\mathbf 1_{\{k\}}(i)
    \equiv-1\pmod3,
\]
so $\mathbf 1_I(i)=1$. This would force $S\subset I$, contrary
to $|S|=3$ and $|I|=2$. Thus the image of $L_\theta$ has no
$\varepsilon_S$-weight space.
Since $L_\theta$ is $A$-equivariant by
Lemma~\ref{lem:Lefschetz-equivariant}, its image has no
$\varepsilon_S$-component. Lemma~\ref{lem:numerator-weight-space}
therefore shows that $q_S$ has a nonzero class in $Q$ and
that $Q_{\varepsilon_S}=\mathbb C[q_S]$.
\end{proof}

In particular, for each of the ten distinct weights in $R_F^3$
we have
\begin{equation}\label{eq:corresponding-weight-lines}
    (R_F^3)_{\varepsilon_S}=\mathbb C v_S,
    \qquad
    Q_{\varepsilon_S}=\mathbb C[q_S].
\end{equation}

%------------------------------------------------
\subsection{The stabilizer}
\label{subsec:stabilizer}
%------------------------------------------------

The construction below follows the stabilizer method for representations
of semidirect products with an abelian normal subgroup, as described
in \cite[\S 8.2]{Serre}. We fix
\[
    S_0:=\{0,1,2\},\qquad v_{S_0}=x_0x_1x_2,
\]
and consider the character $\varepsilon_{S_0}$ of $A$. Since the ten
characters $\varepsilon_S$, $|S|=3$, are distinct, the stabilizer of
$\varepsilon_{S_0}$ in $S_5$ is precisely the subgroup preserving $S_0$:
\begin{equation}\label{eq:stabilizer-S0}
    K:=\operatorname{Stab}_{S_5}(\varepsilon_{S_0})
     =\operatorname{Stab}_{S_5}(S_0)
     \simeq S_3\times S_2.
\end{equation}
The $S_5$-orbit of $v_{S_0}$ consists of the ten basis vectors
$v_S$ of $R_F^3$. Every $k\in K$ fixes $v_{S_0}$, since the action
on $R_F^3$ is the ordinary polynomial action.

By Lemma~\ref{lem:Lefschetz-no-epsilonS}, the corresponding weight
space of the primitive quotient is
\[
    Q_{\varepsilon_{S_0}}=\mathbb C[q_{S_0}],\qquad
    q_{S_0}=(f_0\wedge f_1\wedge f_2)\otimes(e_3\wedge e_4).
\]

\begin{lemma}\label{lem:stabilizer-qS0}
The stabilizer $K$ acts trivially on $\mathbb C[q_{S_0}]$.
\end{lemma}

\begin{proof}
Write $k=(\alpha,\beta)\in S_3\times S_2$, where $\alpha$
permutes $\{0,1,2\}$ and $\beta$ permutes $\{3,4\}$. Reordering
the factors of the two exterior products contributes the sign
\[
    \operatorname{sgn}(\alpha)\operatorname{sgn}(\beta)
    =\operatorname{sgn}(k).
\]
On the other hand, $M$ and $N$ are both twisted by the determinant
character $\chi$. Since $q_{S_0}$ has five factors, the total twist
is $\chi(k)^5=\chi(k)=\operatorname{sgn}(k)$. Thus
\[
    k\cdot q_{S_0}
    =\operatorname{sgn}(k)\chi(k)\,q_{S_0}
    =q_{S_0},
\]
and the same holds for its class in $Q$.
\end{proof}

Consequently, the lines $\mathbb C v_{S_0}$ and
$\mathbb C[q_{S_0}]$ carry the same character of $A\rtimes K$:
$A$ acts on both by $\varepsilon_{S_0}$ and $K$ acts trivially.
This is the compatibility that will allow us to extend
$v_{S_0}\mapsto[q_{S_0}]$ to a $G$-equivariant map.

%------------------------------------------------
\subsection{The equivariant homomorphism and its rank}
\label{subsec:equivariant-homomorphism}
%------------------------------------------------

We now study the space of $G$-equivariant homomorphisms $\operatorname{Hom}_G(R_F^3,Q).$

Fix $S_0=\{0,1,2\}$. For each subset
$S\subset\{0,\ldots,4\}$ with $|S|=3$, choose
$\sigma\in S_5$ such that $\sigma(S_0)=S$,
and consider the linear map
$\phi:R_F^3\to Q$ defined on the basis
$\{v_S:|S|=3\}$ by
\begin{equation}\label{eq:phi-definition}
    \phi(v_S):=\sigma\cdot[q_{S_0}].
\end{equation}
We shall show that this definition is independent
of the choice of $\sigma$ and determine the space
of $G$-equivariant homomorphisms.

\begin{theorem}\label{thm:HomG-Fermat}
The map $\phi$ is well defined, nonzero and
$G$-equivariant. Moreover,
\[
    \operatorname{Hom}_G(R_F^3,Q)=\mathbb C\phi,
    \qquad
    \operatorname{rank}(\phi)=10.
\]
In particular, every nonzero $G$-equivariant
homomorphism from $R_F^3$ to $Q$ is injective.
\end{theorem}

\begin{proof}
We first check that \eqref{eq:phi-definition} is well defined.
If $\sigma(S_0)=\sigma'(S_0)=S$, then
$\sigma'=\sigma k$ for some $k\in K$. Since $K$ fixes
$[q_{S_0}]$ by Lemma~\ref{lem:stabilizer-qS0}, we have
\[
    \sigma'\cdot[q_{S_0}]
    =\sigma k\cdot[q_{S_0}]
    =\sigma\cdot[q_{S_0}].
\]
Thus \eqref{eq:phi-definition} determines a linear map on the basis
$\{v_S:|S|=3\}$, and it is nonzero because
$\phi(v_{S_0})=[q_{S_0}]\ne0$.

The map is $S_5$-equivariant: for $\tau\in S_5$ and
$S=\sigma(S_0)$, the polynomial action gives
$\tau\cdot v_S=(\tau\sigma)\cdot v_{S_0}$, then
\[
    \phi(\tau\cdot v_S)
    =\tau\sigma\cdot[q_{S_0}]
    =\tau\cdot\phi(v_S).
\]
Since $[q_{S_0}]$ has $A$-weight
$\varepsilon_{S_0}$ and $\sigma(S_0)=S$,
the vector
$\phi(v_S)=\sigma\cdot[q_{S_0}]$
has $A$-weight $\varepsilon_S$.
Therefore $\phi(v_S)\in Q_{\varepsilon_S}.$
As $v_S$ has the same
$A$-weight, for every $a\in A$ we obtain
\[
    \phi(a\cdot v_S)
    =\varepsilon_S(a)\phi(v_S)
    =a\cdot\phi(v_S).
\]
Hence $\phi$ is $G=A\rtimes S_5$-equivariant.

Now let $\psi\in\operatorname{Hom}_G(R_F^3,Q)$.
By $A$-equivariance and
Lemma~\ref{lem:Lefschetz-no-epsilonS},
\[
    \psi(v_{S_0})\in Q_{\varepsilon_{S_0}}
    =\mathbb C[q_{S_0}].
\]
Thus $\psi(v_{S_0})=c\phi(v_{S_0})$ for some
$c\in\mathbb C$. Since the vectors $v_S=\sigma\cdot v_{S_0}$
form a basis and both maps are $S_5$-equivariant, it follows that
$\psi(v_S)=c\phi(v_S)$ for every $S$. Therefore
$\psi=c\phi$.

It remains to prove that $\phi$ is injective.
For each $S$, the vector $\phi(v_S)$ is nonzero
and lies to the weight space $Q_{\varepsilon_S}$.
Since the characters $\varepsilon_S$ are pairwise
distinct, the images of the basis vectors $v_S$
are linearly independent. Since $\dim R_F^3=10$, we conclude that
$\operatorname{rank}(\phi)=10$.
\end{proof}

%------------------------------------------------
\section{Nonvanishing and maximal rank}
\label{sec:nonvanishing}
%------------------------------------------------

Theorem~\ref{thm:HomG-Fermat} and
\eqref{eq:canonical-HomG-Fermat} show that $\overline\nabla_0\nu_F=c\phi$ for some $c \in \mathbb C$.
To prove maximal rank at the Fermat cubic, it remains to show
that $c\ne0$. By \eqref{eq:delta-nonzero-canonical}, it is enough to prove the nonvanishing of the infinitesimal invariant of $\nu_F$ at the Fermat point.

We recall the criterion of Collino--Naranjo--Pirola
\cite[Proposition~5.9(c)]{CNP}. For a smooth cubic threefold
$V\subset\mathbb P^4$ and $0\ne\xi\in R^3(V)$, set
\[
    K_1(\xi):=\ker\bigl(R^1(V)
          \xrightarrow{\,\cdot\xi\,}R^4(V)\bigr).
\]
Let $\mathbb G:=\operatorname{Grass}(2,5)$
be the Grassmannian of lines in \(\mathbb P^4\). For
\(r\in\mathbb G\), denote by $\ell_r\subset\mathbb P^4$
the corresponding line and set
\[
    W(r)
    :=
    \left\{
        L\in R^1(V)
        \mid
        L|_{\ell_r}=0
    \right\}.
\]
Thus \(W(r)\) is the three-dimensional vector space of hyperplanes
containing \(\ell_r\).
A line $\ell_r$ is called \emph{special} if
$K_1(\xi)=W(r)$ for some $0\ne\xi\in R^3(V)$; such a
class $[\xi]$ is unique by \cite[Lemma~3.4]{CNP}.

\begin{proposition}[Collino--Naranjo--Pirola
    {\cite[Proposition~5.9(c)]{CNP}}]
\label{prop:CNP-nonvanishing}
Suppose $K_1(\xi)=W(r)$ and that $\ell_r$ is not contained
in $V$ but intersects $V$ in three distinct points. Then
\[
    \delta_V(\nu_F)\ne0.
\]
\end{proposition}

Applying Proposition \ref{prop:CNP-nonvanishing} to the Fermat cubic threefold, we have the final result.

\begin{theorem}\label{thm:fermat-nonvanishing}
At the Fermat cubic $V_F$, the infinitesimal invariant of the
Fano normal function is nonzero. Consequently,
\[
    \operatorname{rk}_{0}(\nu_F)
    =\operatorname{rank}_{\mathbb C}
      (\overline\nabla_0\nu_F)=10.
\]
In particular, $\nu_F$ has maximal rank on a dense open subset
of the moduli space of smooth cubic threefolds.
\end{theorem}

\begin{proof}
Consider the line
\[
    \ell_r:=\{x_2=x_3=x_4=0\}\subset\mathbb P^4,
    \qquad W(r)=\langle x_2,x_3,x_4\rangle,
\]
and the nonzero deformation class
$\xi:=x_2x_3x_4\in R_F^3$. For
$L=\sum_{i=0}^{4}a_ix_i\in R_F^1$, the relations
$x_i^2=0$ in $R_F^\bullet$ give
\[
    L\xi
    =a_0x_0x_2x_3x_4+a_1x_1x_2x_3x_4.
\]
The two monomials on the right are linearly independent in
$R_F^4$. Hence
\[
    K_1(\xi)=\langle x_2,x_3,x_4\rangle=W(r),
\]
so $\ell_r$ is special. On this line the equation of $V_F$
restricts to
\[
    E_F|_{\ell_r}=x_0^3+x_1^3.
\]
This cubic is not identically zero and has three distinct zeros
in $\ell_r\simeq\mathbb P^1$. Thus $\ell_r\not\subset V_F$
and $\ell_r\cap V_F$ consists of three distinct points.
Proposition~\ref{prop:CNP-nonvanishing} now gives
$\delta_0(\nu_F)\ne0$.

By \eqref{eq:delta-nonzero-canonical},
$\overline\nabla_0\nu_F\ne0$. Since 
$\overline\nabla_0\nu_F \in \operatorname{Hom}_G(R_F^3,Q)=\mathbb C\phi$,
Theorem~\ref{thm:HomG-Fermat} implies that it has rank $10$.

Finally, by Remark~\ref{remark dense open},
maximal rank at the Fermat cubic implies maximal
rank on a dense open subset of the moduli space.
\end{proof}

\begin{remark}
The preceding argument also gives a more precise
description of the canonical derivative at the Fermat
cubic. As recalled in Section~2.3, the infinitesimal
invariant determines the canonical derivative only
modulo the image of the infinitesimal variation of
Hodge structure.

At the Fermat cubic, however, this ambiguity
disappears among $G$-equivariant homomorphisms.
Indeed, by Theorem~\ref{thm:HomG-Fermat},
\[
    \operatorname{Hom}_G(R_F^3,Q)=\mathbb C\varphi.
\]
Since $\delta_0(\nu_F)\neq0$, the class of
$\varphi$ modulo the image of the infinitesimal
variation of Hodge structure is nonzero.
Consequently, the infinitesimal invariant
$\delta_0(\nu_F)$ uniquely determines
$\overline\nabla_0\nu_F$ among
$G$-equivariant homomorphisms from $R_F^3$ to $Q$.
\end{remark}

\subsection*{Acknowledgements}

I am grateful to Gian Pietro Pirola for suggesting this problem
during my PhD studies at the University of Pavia.
I would also like to thank Elisabetta Colombo, Paola Frediani,
and Juan Carlos Naranjo for stimulating discussions on this topic. The author is a member of GNSAGA--INdAM.

\end{document}